\documentclass[dvipdfmx]{amsart}
\usepackage{tikz}
\usetikzlibrary{arrows.meta}
\usepackage{tkz-graph}
\usetikzlibrary{arrows}

\usepackage{amssymb}
\usepackage{amsmath}
\usepackage{amsthm}
\usepackage{times}
\usepackage{graphicx}

\usepackage[curve,all]{xy}
 \xyoption{curve}
 \xyoption{line}
\xyoption{arc}

\newtheorem{theorem}{Theorem}[section]
\newtheorem{lemma}[theorem]{Lemma}
\newtheorem{proposition}[theorem]{Proposition}
\theoremstyle{definition}
\newtheorem{definition}[theorem]{Definition}
\newtheorem{example}[theorem]{Example}
\theoremstyle{remark}
\newtheorem{remark}[theorem]{Remark}

\numberwithin{equation}{section}

\def\PtAv{v}
\def\PtAw{w}
\def\hVar{z}

\def\Ker{{\text{Ker}}}
\def\deg{{\text{deg}}}

\begin{document}

\title[Counting Richelot isogenies between superspecial abelian surfaces]{Counting Richelot isogenies between \\ superspecial abelian surfaces}

\author{Toshiyuki Katsura}
\address{Graduate School of Mathematical Sciences, The University of Tokyo, 
Meguro-ku, Tokyo 153-8914, Japan}
\email{tkatsura@ms.u-tokyo.ac.jp}

\author{Katsuyuki Takashima}
\address{Information Technology R\&D Center, Mitsubishi Electric, Kamakura-shi, Kanagawa 247-8501, Japan}
\email{Takashima.Katsuyuki@aj.MitsubishiElectric.co.jp}

\begin{abstract}
Castryck, Decru, and Smith used superspecial genus-2 curves and their Richelot isogeny graph for basing genus-2 isogeny cryptography, and recently, Costello and Smith devised an improved isogeny path-finding algorithm in the genus-2 setting. In order to establish a firm ground for the cryptographic construction and analysis, we give a new characterization of {\em decomposed Richelot isogenies} in terms of {\em involutive reduced automorphisms} of genus-2 curves over a finite field, and explicitly count such decomposed (and non-decomposed) Richelot isogenies between {\em superspecial} principally polarized abelian surfaces. As a corollary, we give another algebraic geometric proof of Theorem 2 in the paper of Castryck et al.
\end{abstract}



\maketitle

\section{Introduction}

Isogenies of supersingular elliptic curves are widely studied as one candidate for post-quantum cryptography, e.g., \cite{CLG09,FJP14,SIKE20,CLMPR18}. Recently, several authors have extended the cryptosystems to higher genus isogenies, especially the genus-2 case \cite{Tak17,FT19,CDS,CS20}.

Castryck, Decru, and Smith \cite{CDS} showed that {\em superspecial} genus-2 curves and their isogeny graphs give a correct foundation for constructing genus-2 isogeny cryptography. 
The recent cryptanalysis by Costello and Smith \cite{CS20} employed the subgraph whose vertices consist of decomposed principally polarized abelian varieties, hence it is important to study the subgraph in cryptography.

Castryck et al.\,also presented concrete algebraic formulas for computing $(2,2)$-isogenies by using the Richelot construction. In the genus-2 case, the isogenies may have decomposed principally polarized abelian surfaces 
as codomain, and we call them decomposed isogenies. 
In \cite{CDS}, the authors gave explicit formulas for the decomposed isogenies and a theorem stating that the number of decomposed Richelot isogenies outgoing from the Jacobian $J(C)$ of a superspecial curve $C$ of genus 2 is {\em at most six} (Theorem 2 in \cite{CDS}), but they {\em do not precisely determine} this number. Moreover, their proof is {\em computer-aided}, that is, using the Gr$\ddot{\rm o}$bner basis computation.     

Therefore, we revisit the isogeny counting based on an intrinsic algebraic geometric characterization. 
In 1960, Igusa \cite{I} classified the curves of genus 2 with given
reduced groups of automorphisms, and in 1986, Ibukiyama, Katsura, and Oort \cite{IKO} explicitly counted such superspecial curves according to the classification.  
Based on the classical results, we first count the number of Richelot isogenies from a superspecial Jacobian to decomposed surfaces (Cases (0)--(6) in Section \ref{Sec:NumberOfLongElements}) in terms of {\em involutive (i.e., of order 2) reduced automorphisms} which are called long elements. As a corollary, we give an algebraic geometric proof of Theorem 2 in \cite{CDS} together with a {\em precise count of 
decomposed Richelot isogenies} (Remark \ref{Rem:Thm3CDS}).
Moreover, by extending the method, we also count the total number of (decomposed) Richelot isogenies up to isomorphism outgoing from irreducible superspecial curves of genus 2 (resp.\,decomposed principally polarized superspecial abelian surfaces) in Theorem \ref{Thm:RichelotFromNonDecomp} (resp.\,Theorem \ref{Thm:RichelotFromDecomp}).

Our paper is organized as follows: Section \ref{Sec:Preliminaries} gives mathematical preliminaries including the Igusa classification and the Ibukiyama--Katsura--Oort curve counting. Section \ref{Sec:RichelotIsogeny} presents an abstract description of Richelot isogenies and Section \ref{Sec:DecomposedRichelot} gives the main characterization of decomposed Richelot isogenies in terms of reduced groups of 
automorphisms. Section \ref{Sec:NumberOfLongElements} counts the number of long elements of order 2 in reduced groups of automorphisms 
based on the results in Section \ref{Sec:DecomposedRichelot}. 
Section \ref{Sec:Counting} gives the total numbers of (decomposed) Richelot isogenies outgoing from the irreducible superspecial curves of genus 2 and products of two elliptic curves, respectively.  
Section \ref{Sec:Examples} gives some examples in small characteristic. Finally, Section \ref{Sec:ConcludingRemarks} gives a concluding remark. 

\vspace*{0.1cm}
We use the following notation:
For an abelian surface $A$,
$A[n]$ denotes the group of $n$-torsion points of $A$,
$A^t$ the dual of $A$, 
${\rm NS}(A)$ the N\'eron--Severi group of $A$, and
$T_\PtAv$ the translation by an element $\PtAv$ of $A$.
For a nonsingular projective variety $X$,
$D \sim D'$ (resp.\,$D \approx D'$) denotes linear equivalence (resp.\,numerical equivalence) for divisors $D$ and $D'$ on $X$, and
$id_X$ the identity morphism of $X$.

\newcounter{num}
\setcounter{num}{14}
\vspace*{0.1cm}
\noindent
{\bf Acknowledgements.} The authors would like to thank anonymous reviewers of ANTS-\Roman{num} 
for their careful reading and useful suggestions for revising our paper, and E.~Florit and B.~Smith for their useful comments, in particular, for the correction of the figure in the case of $p = 11$ in Subsection \ref{Sec:Examples11}.
Research of the first author is partially supported by JSPS Grant-in-Aid for Scientific Research (C) No.\,20K03530.
Research of the second author is partially supported by JST CREST Grant Number  JPMJCR14D6, Japan.


\section{Preliminaries}
\label{Sec:Preliminaries}

Let $k$ be an algebraically closed field of characteristic $p > 5$. An abelian 
surface $A$ defined over $k$
is said to be superspecial if $A$ is isomorphic to $E_{1}\times E_{2}$ 
with $E_{i}$ supersingular elliptic curves $(i = 1, 2)$. 
Since for any supersingular
elliptic curves $E_{i}$ $(i = 1, 2, 3, 4)$ we have an isomorphism 
$E_{1}\times E_{2} \cong E_{3}\times E_{4}$  
(cf.~Shioda \cite[Theorem 3.5]{S}, for instance), this notion does not depend
on the choice of supersingular elliptic curves. For a nonsingular projective curve
$C$ of genus 2, we denote by 
$(J(C),C)$ the canonically polarized Jacobian variety of $C$. 
The curve $C$ is said to be superspecial if $J(C)$ is superspecial as an abelian surface.
We denote by ${\rm Aut} (C)$
the group of automorphisms of $C$. Since $C$ is hyperelliptic, $C$ has 
the hyperelliptic involution $\iota$ such that 
the quotient curve $C/\langle \iota \rangle$
is isomorphic to the projective line ${\bf P}^{1}$:
$$
         \psi : C \longrightarrow {\bf P}^{1}.
$$
There exist 6 ramification points on $C$. We denote them by $P_{i}$
$(1 \leq i \leq 6)$. Then, $Q_{i} = \psi(P_{i})$'s are the branch points of $\psi$
on ${\bf P}^{1}$. 
The group $\langle \iota \rangle$
is a normal subgroup of ${\rm Aut} (C)$. We put ${\rm RA}(C) \cong {\rm Aut} (C)/\langle \iota\rangle$ and we call it the reduced group of automorphisms of $C$. 
We call an element of ${\rm RA}(C)$ a reduced automorphism of $C$.
For $\sigma \in {\rm RA}(C)$, $\tilde{\sigma}$ is an element of ${\rm Aut}(C)$ such that 
$\tilde{\sigma} ~{\rm mod}~ \langle \iota \rangle = \sigma$.
\begin{definition}
An element $\sigma \in {\rm RA}(C)$ of order 2 is said to be long
if $\tilde{\sigma}$ is of order 2. Otherwise, 
an element $\sigma \in {\rm RA}(C)$ of order 2 is said to be short
(cf.~Katsura--Oort \cite[Definition 7.15]{KO2}).
\end{definition}
This definition does not depend on the choice of $\tilde{\sigma}$.

\begin{lemma}\label{freely}
If an element $\sigma \in {\rm RA}(C)$ of order 2 acts freely on 6 branch points,
then $\sigma$ is long.
\end{lemma}

\proof{By a suitable choice of coordinate $x$ of ${\bf A}^{1}\subset {\bf P}^{1}$, 
taking 0 as a fixed point of $\sigma$, 
we may assume $\sigma (x) = -x$, and $Q_1 = 1$, $Q_2 = -1$, $Q_3 = a$, $Q_4 = -a$,
$Q_5 = b$, $Q_6 = -b$ $(a\neq 0, \pm1; b\neq 0, \pm1; a\neq \pm b)$. Then, the curve is defined by
$$
     y^2 = (x^2 - 1)(x^2 - a^2)(x^2 - b^2),
$$
and $\tilde{\sigma}$ is given by
$
      x\mapsto -x, y \mapsto \pm y.
$
Therefore, $\tilde{\sigma}$ is of order 2.
\qed}

\begin{lemma}\label{replace}
If ${\rm RA}(C)$ has an element $\sigma$ of order 2, then there exists
a long element $\tau \in {\rm RA}(C)$ of order 2.
\end{lemma}

\proof{If $\sigma$ acts freely on 6 branch points, then by Lemma \ref{freely},
$\sigma$ itself is a long element of order 2. We assume that 
the branch point $Q_{1} = \psi(P_1)$
is a fixed point of $\sigma$. Since $\sigma$ is of order 2, it must have
one more fixed point among the branch points, say $Q_{2} = \psi(P_2)$.
By a suitable choice of coordinate $x$ of ${\bf A}^{1}\subset {\bf P}^{1}$, 
we may assume $Q_{1} = 0$ and $Q_{2} = \infty$. We may also assume $Q_{3} = 1$.
Then, $\sigma$ is given by $x \mapsto -x$ and
the six branch points are $0$, $1$, $-1$, $a$, $-a$, $\infty$ ($a \neq \pm 1$).
The curve $C$ is given by
$$
     y^2 = x(x^2 - 1)(x^2 - a^2) \quad (a\neq 0, \pm1).
$$
We consider an element $\tau \in {\rm Aut}({\bf P}^{1})$ defined by
$ 
       x \mapsto \frac{a}{x}
$.
Then, we have an automorhisms $\tilde{\tau}$ of $C$ defined by
$
         x \mapsto  \frac{a}{x}, y \mapsto \frac{a\sqrt{a}y}{x^3}
$.
Therefore, we see $\tau \in {\rm RA}(C)$. Since $\tilde{\tau}$ is of order 2,
$\tau$ is long.
\qed}
\vspace*{0.1cm}

${\rm RA}(C)$ acts
on the projective line ${\bf P}^{1}$ as a subgroup of ${\rm PGL}_{2}(k)$.
The structure of ${\rm RA}(C)$ is classified as follows 
(cf.~Igusa \cite[p. 644]{I}, and Ibukiyama--Katsura--Oort \cite[p. 130]{IKO}):
$$
(0)\,0, ~(1)\, {\bf Z}/2{\bf Z},~(2)\, S_{3},~(3)\, {\bf Z}/2{\bf Z}\times {\bf Z}/2{\bf Z},
~(4)\, D_{12}, ~ (5)\, S_{4},~(6)\, {\bf Z}/5{\bf Z}.
$$
We denote by $n_{i}$ the number of superspecial curves of genus 2 whose reduced
group of automorphisms is isomorphic to the group $(i)$.
Then, $n_{i}$'s are given as follows (cf.~Ibukiyama--Katsura--Oort 
\cite[Theorem 3.3]{IKO}):
\begin{enumerate}
\item[$(0)$] $n_{0}= (p -1)(p^2 - 35p + 346)/2880  -\{1 - (\frac{-1}{p})\}/32 
-\{1 - (\frac{-2}{p})\}/8$ $- \{1 - (\frac{-3}{p})\}/9$
$+
\left\{
\begin{array}{ll}
0 & {\rm if}~p \equiv 1, 2 ~{\rm or}~3 ~ ({\rm mod}~ 5),\\
-1/5 & {\rm if}~p \equiv 4 ~ ({\rm mod}~ 5),
\end{array}
\right.
$
\item[$(1)$] $n_{1}= (p-1)(p-17)/48+ \{1 - (\frac{-1}{p})\}/8  
                       +\{1 - (\frac{-2}{p})\}/2 + \{1 - (\frac{-3}{p})\}/2$,
\item[$(2)$] $n_{2}= (p - 1)/6 -\{1 - (\frac{-2}{p})\}/2 - \{1 - (\frac{-3}{p})\}/3$,
\item[$(3)$] $n_{3}= (p - 1)/8 -\{1 - (\frac{-1}{p})\}/8 - 
                       \{1 - (\frac{-2}{p})\}/4 - \{1 - (\frac{-3}{p})\}/2$,
\item[$(4)$] $n_{4}= \{1 - (\frac{-3}{p})\}/2$,
\item[$(5)$] $n_{5}= \{1 - (\frac{-2}{p})\}/2$,
\item[$(6)$] $n_{6}=$
$
\left\{
\begin{array}{ll}
0 & {\rm if}~p \equiv 1, 2 ~{\rm or}~3 ~ ({\rm mod}~ 5),\\
1 & {\rm if}~p \equiv 4 ~ ({\rm mod}~ 5).
\end{array}
\right.
$
\end{enumerate}

Here, for a prime number $q$ and an integer $a$, $(\frac{a}{q})$
is the Legendre symbol. The total number $n$ of superspecial curves of genus 2 is
given by
$$
\begin{array}{cl}
n &= n_{0} + n_{1} + n_{2} + n_{3} + n_{4} + n_{5} + n_{6}\\
  &= (p - 1)(p^2 + 25p + 166)/2880 - \{1 - (\frac{-1}{p})\}/32  
                       +\{1 - (\frac{-2}{p})\}/8 \\
  &\quad + \{1 - (\frac{-3}{p})\}/18+\left\{
\begin{array}{ll}
0 & {\rm if}~p \equiv 1, 2 ~{\rm or}~3 ~ ({\rm mod}~ 5),\\
4/5 & {\rm if}~p \equiv 4 ~ ({\rm mod}~ 5).
\end{array}
\right.
\end{array}
$$

For an abelian surface $A$, we have $A^{t} = {\rm Pic}^{0}(A)$ (Picard variety of $A$), and
for a divisor $D$ on $A$,
there exists a homomorphism
$$
\begin{array}{cccc}
   \varphi_{D}  :& A   &  \longrightarrow &  A^{t} \\
      &   \PtAv   &  \mapsto  & T_{\PtAv}^{*}D - D.
\end{array}
$$
If $D$ is ample, then $\varphi_{D}$ is surjective, i.e., an isogeny.
We know $(D\cdot D)^2 = 4~\deg~\varphi_{D}$.
We set $K(D) = \Ker ~\varphi_{D}$. If $D$ is ample, then $K(D)$
is finite and there is a non-degenerate alternating bilinear
form $e^{D}(\PtAv, \PtAw)$ on $K(D)$ (cf.~Mumford \cite[Section 23]{M}). 
Let $G$ be an isotropic subgroup scheme of $K(D)$
with respect to $e^{D}(\PtAv, \PtAw)$. In case $D$ is ample, $G$ is finite and 
we have an isogeny
$$
   \pi : A \longrightarrow A/G.
$$

The following theorem is due to Mumford 
\cite[Section 23, Theorem 2, Corollary]{M}:

\begin{theorem}\label{descent}
Let $G$ be an isotropic subgroup scheme of $K(D)$. Then, there exists
a divisor $D'$ on $A/G$ such that $\pi^{*}D' \sim D$. 
\end{theorem}

Let $n$ be a positive integer which is prime to $p$.
Then, we have the Weil pairing 
$e_{n} : A[n]\times A^{t}[n] \longrightarrow \mu_{n}$. 
Here, $\mu_{n}$ is the multiplicative group of order $n$.
By Mumford \cite[Section 23 ``Functional Properties of $e^{L}$ (5)'']{M},
we have the following.

\begin{lemma}
For $\PtAv \in A[n]$ and $\PtAw \in \varphi_{D}^{-1}(A^{t}[n])$, we have
$$
       e_{n}(\PtAv, \varphi_{D}(\PtAw)) = e^{nD}(\PtAv, \PtAw).
$$
\end{lemma}
If $D$ is a principal polarization, the homomorphism 
$\varphi_{D} : A \longrightarrow A^{t}$ is an isomorphism.
Therefore, by this identification we can identify the pairing $e^{nD}$
with the Weil pairing $e_{n}$.


\section{Richelot isogenies}
\label{Sec:RichelotIsogeny}

We recall the abstract description of Richelot isogenies.
(For the concrete construction of Richelot isogenies, see 
Smith \cite{Sm}, Castryck--Decru--Smith \cite[Section 3]{CDS}, 
for instance.)

Let $A$ be an abelian surface with a principal polarization $C$.
Then, we may assume that $C$ is effective, and we have 
the self-intersection number $C^2 = 2$.
It is easy to show (or as was shown by A. Weil) that there are two cases
for effective divisors with self-intersection 2 on an abelian surface $A$:

(1) There exists a nonsingular curve $C$ of genus 2 such that $A$ is 
isomorphic to the Jacobian variety $J(C)$ of $C$ and that $C$ is the divisor 
with self-intersection 2.
In this case, $(J(C), C)$ is said to be non-decomposed. 

(2) There exist two elliptic curves $E_{1}$, $E_{2}$ with $(E_{1}\cdot E_{2}) = 1$
such that $E_1\times \{0\} + \{0\} \times E_2$ is a divisor with self-intersection 2 and that $A \cong E_{1}\times E_{2}$.
In this case, $(A, E_1\times \{0\} + \{0\} \times E_2)$ is said to be decomposed.

Since $\varphi_{C}$ is an isomorphism by the fact that $C$ is a principal
polarization, we have 
$K(2C) = \Ker~\varphi_{2C} = \Ker~2\varphi_{C}= A[2]$.
Let $G$ be a maximal isotropic subgroup of $K(2C)= A[2]$ with respect to 
the pairing $e^{2C}$. Since we have $|G|^{2} = |A[2]| = 2^4$  
(cf.~Mumford \cite[Section 23, Theorem 4]{M}), 
we have $|G| = 4$
and $G \cong {\bf Z}/2{\bf Z} \times {\bf Z}/2{\bf Z}$. We have a quotient homomorphism
$$
      \pi : A \longrightarrow A/G.
$$
By Theorem \ref{descent}, there exists a divisor $C'$ on $A/G$ 
such that 
$2C \sim \pi^{*}C'$.
Since $\pi$ is a finite morphism and $2C$ is ample, 
we see that $C'$ is also ample.
We have the self-intersection number $(2C \cdot 2C) = 8$, and we have
$$
8 = (2C \cdot 2C) = (\pi^{*}C'\cdot \pi^{*}C') = \deg~\pi~(C'\cdot C') 
= 4(C'\cdot C').
$$
Therefore, we have $(C'\cdot C') = 2$, that is, $C'$ is a principal polarization
on $A/G$. By the Riemann--Roch theorem of an abelian surface for ample divisors,
we have 
$$
    \dim {\rm H}^{0}(A/G, {\mathcal O}_{A/G}(C')) = (C'\cdot C')/2 = 1.
$$
Therefore, we may assume $C'$ is an effective divisor. 

Using  these facts, we see that $C'$ is either a nonsingular curve of genus 2
or $E_{1} \cup E_{2}$ with elliptic curves $E_{i}$ ($i = 1, 2$) which
intersect each other transversely.
In this situation, the correspondence from $(A, C)$ to $(A/G, C')$
is called a Richelot isogeny. 
We consider a triple $(A, C, G)$ with maximal isotropic subgroup $G \subset A[2]$ with respect to the pairing $e^{2C}$, and the corresponding Richelot isogeny $\pi$ from $(A, C, G)$ to $(A/G, C', G')$ with maximal isotropic subgroup $G'= \pi(A[2])$.
Then, it is easy to see that for the Richelot isogeny 
$\pi' : (A/G, C') \longrightarrow ((A/G)/G', C'')$, the principally polarized abelian surface $((A/G)/G', C'', G'')$ with maximal isotropic subgroup $G''= \pi'((A/G)[2])$ is isomorphic to the original $(A, C, G)$.

Now, we consider the case where $A$ is a superspecial abelian surface. 
Then,
since $\pi$ is separable, 
$A/G$ is also a superspecial abelian surface. We will use this fact freely.

From here on, for abelian surface $E_1 \times E_2$ with elliptic curves $E_i$
($i = 1, 2$) we donote by $E_1 + E_2$ 
the divisor $E_1\times \{0\} + \{0\} \times E_2$, if no confusion occurs.
We sometimes call $E_1 \times E_2$ a principally polarized abelian surface.
In this case, the principal polarization on $E_1 \times E_2$ is given by 
$E_1 + E_2$.

\begin{definition}
\label{Def:IsomorphismRichelot}
Let $(A, C)$, $(A', C')$ and $(A'', C'')$ be principally polarized
abelian surfaces with principal polarizations $C$, $C'$, $C''$, respectively.
The Richelot isogeny $\pi : A \longrightarrow A'$
is said to be isomorphic to the Richelot isogeny
$\varpi : A \longrightarrow A''$ if there exist an
automorphism $\sigma \in A$ with $\sigma^{*}C \approx C$ and an isomorphism
$g : A' \longrightarrow A''$ with $g^{*}C'' \approx C'$ 
such that the
following diagram commutes:
$$
\begin{array}{rrcl}
  \sigma :&    A ~& \longrightarrow & A \\
     &  \pi \downarrow ~ &       & ~\downarrow \varpi \\
  g  : &  A' & \longrightarrow   & A''.
\end{array}
$$
\end{definition}


\section{Decomposed Richelot isogenies}
\label{Sec:DecomposedRichelot}

In this section, we use the same notation as in Section 3.

\begin{definition} Let $A$ and $A'$ be abelian surfaces 
with principal polarizations
$C$, $C'$, respectively. A Richelot isogeny 
$A \longrightarrow A'$ is said to be decomposed 
if $C'$ consists of two elliptic curves. Otherwise, the Richelot isogeny
is said to be non-decomposed.
\end{definition}

\begin{example}\label{Igusa}
Let $C_{a, b}$ be a nonsingular projective model of the curve of genus 2
defined by the equation
$$
 y^{2} = (x^{2} -1)(x^{2} - a)(x^{2} - b)\quad (a \neq 0, 1; b \neq 0, 1; a\neq  b).
$$
Let $\iota$ be the hyperelliptic involution defined by $x \mapsto x,~y \mapsto -y$.
${\rm RA}(C_{a, b})$ has an element of order 2 defined by 
$$
\sigma : x \mapsto -x, y \mapsto y.
$$
We put $\tau = \iota \circ \sigma$. We have two elliptic curves $E_{\sigma}= C_{a, b}/\langle \sigma \rangle$ and $E_{\tau}= C_{a, b}/\langle \tau \rangle$.
The elliptic curve $E_{\sigma}$ is isomorphic to an elliptic curve
$
   E_{\lambda} :  y^2 = x(x - 1)(x - \lambda)
$
with
\begin{equation}\label{equation1}
  \lambda = (b - a)/(1 - a)
\end{equation}
and the elliptic curve $E_{\tau}$ is isomorphic to an elliptic curve
$
    E_{\mu} :  y^2 = x(x - 1)(x - \mu)
$
with
\begin{equation}\label{equation2}
 \mu = (b - a)/b(1 - a).
\end{equation}
The map given by (\ref{equation1}) and (\ref{equation2}) yields a bijection
$$
\begin{array}{c}
\{(a, b) \mid a, b\in k; a \neq0, 1; b \neq 0, 1; a\neq b, {\rm and} ~J(C_{a, b})~{\rm is ~superspecial}\} \\
\longrightarrow \{(\lambda, \mu) \mid \lambda, \mu \in k; \lambda \neq \mu ; E_{\lambda}, E_{\mu} ~{\rm are ~supersingular}\}
\end{array}
$$
(for the details, see Katsura--Oort \cite[p. 259]{KO}).
We have a natural morphism $C_{a, b} \longrightarrow  E_{\sigma} \times E_{\tau}$
and
this morphism induces an isogeny
$$
\pi : J(C_{a, b}) \longrightarrow E_{\sigma} \times E_{\tau}.
$$
By Igusa \cite[p. 648]{I}, we know 
$\Ker~\pi \cong {\bf Z}/2{\bf Z} \times {\bf Z}/2{\bf Z}$ and $\Ker~\pi$
consists of $P_{1} - \sigma (P_{1})$, $P_{3} - \sigma (P_{3})$, $P_{5} - \sigma (P_{5})$ and the zero point. 
Here, $P_{1} =(1, 0)$, $P_{3} =(a, 0)$, $P_{5} =(b, 0)$.
Since $P_{i} - \sigma (P_{i})$ is a divisor 
of order 2, we have $P_{i} - \sigma (P_{i}) \sim \sigma (P_{i}) - P_{i}$.

Comparing the calculation in Castryck--Decru--Smith \cite[Proposition 1 (2)]{CDS} 
with the one in Katsura--Oort \cite[Lemma 2.4]{KO}, we see that
$\pi : J(C_{a, b}) \longrightarrow E_{\sigma} \times E_{\tau}$ is a decomposed Richelot 
isogeny with $C_{a, b}' = E_{\sigma} + E_{\tau}$ (also see 
Katsura--Oort \cite[Proof of Proposition 7.18 (iii)]{KO2}). 
We will use the bijection above to calculate
decomposed Richelot isogenies.
\end{example}

\begin{proposition}\label{main}
Let $C$ be a nonsingular projective curve of genus 2.
Then, the following three conditions are equivalent.

(i) $C$ has a decomposed Richelot isogeny outgoing from $J(C)$.

(ii) ${\rm RA}(C)$ has an element of order 2.

(iii) ${\rm RA}(C)$ has a long element of order 2.
\end{proposition}

\proof{(i) $\Rightarrow$ (ii). By assumption, we have a Richelot isogeny
\begin{equation}\label{Richelot}
     \pi : J(C) \longrightarrow J(C)/G
\end{equation}
such that $G$ is an isotropic subgroup of $J(C)[2]$ with respect to $2C$, and that
$C'$ is a principal polarization consisting of two elliptic curves 
$E_{i}$ ($i = 1, 2$) on $J(C)/G$ with $2C \sim \pi^{*}(E_{1} + E_{2})$. 
Since $C$ is a principal polarization, we have 
an isomorphism $\varphi_C : J(C) \cong  J(C)^{t}$. 
In a similar way, we have $J(C)/G \cong (J(C)/G)^{t}$.
Dualizing (\ref{Richelot}), we have 
$$
       \eta = \pi^{t} : J(C)/G \longrightarrow J(C)
$$
with $J(C)/G \cong E_{1}\times E_{2}$, $C' = E_{1} + E_{2}$ and $\eta^{*}(C) \sim 2(E_{1} + E_{2})$. 
The kernel $\Ker~ \eta$ is an isotropic subgroup of $(E_{1}\times E_{2})[2]$ 
with respect to the divisor $2(E_{1} + E_{2})$. 

Denoting by $\iota_{E_{1}}$ the inversion of $E_{1}$, we set
$$
 \bar{\tau} = \iota_{E_{1}}\times id_{E_{2}}.
$$
Then, $\bar{\tau}$ is an automorphism of order 2 which is not the inversion of 
$E_{1}\times E_{2}$. By the definition, we have
$$
\bar{\tau}^{*}(E_{1} + E_{2}) = E_{1} + E_{2}.
$$
Moreover, since $\Ker~ \eta$ consists of elements of order 2 and 
$\bar{\tau}$ fixes the elements of order 2, 
$\bar{\tau}$ preserves $\Ker~ \eta$. 
Therefore, $\bar{\tau}$ induces an automorphism 
$\tau$ of $J(C) \cong (J(C)/G)/\Ker~ \eta \cong (E_{1}\times E_{2})/\Ker~ \eta$.
Therefore, we have the following diagram:
$$
\begin{array}{ccc}
 E_{1}\times E_{2} & \stackrel{\bar{\tau}}{\longrightarrow}  & E_{1}\times E_{2}\\
          \eta \downarrow  &        &  \downarrow \eta \\
          J(C) & \stackrel{\tau}{\longrightarrow}& J(C).
\end{array}
$$
We have
$$
  \eta^{*}\tau^{*}C= \bar{\tau}^{*}\eta^{*}C \sim \bar{\tau}^{*}(2(E_{1} + E_{2})) = 2(E_{1} + E_{2}).
$$
On the other hand, we have
$$
   \eta^{*}C \sim 2(E_{1} + E_{2}).
$$
Since $\eta^{*}$ is an injective homomorphism from ${\rm NS}(J(C))$ to
${\rm NS}(E_{1}\times E_{2})$, we have $C \approx \tau^{*}C$. 
Therefore, $\tau^{*}C - C$ is an element of 
${\rm Pic}^{0}(J(C)) = J(C)^{t}$. Since $C$ is ample, the homomorphism
$$
\begin{array}{cccc}
   \varphi_{C} : &  J(C) &\longrightarrow & J(C)^{t}\\
     & \PtAv & \mapsto & T_{\PtAv}^{*}C - C
\end{array}
$$
is surjective. 
Therefore, there exists an element $\PtAv \in J(C)$ such that
$$
        T_{\PtAv}^{*}C - C \sim \tau^{*}C - C,
$$
that is, $T_{\PtAv}^{*}C \sim \tau^{*}C$. Since $T_{\PtAv}^{*}C$ is 
a principal polarization, we see 
$$
\dim {\rm H}^{0}(J(C), {\mathcal O}_{J(C)}(T_{\PtAv}^{*}C)) = 1.
$$
Therefore, we have $T_{\PtAv}^{*}C = \tau^{*}C$, that is, 
$T_{-\PtAv}^{*}\tau^{*}C = C$.  Since $\tau$ is of order 2,
we have $(\tau \circ T_{-\PtAv})^2 = T_{-\PtAv -\tau(\PtAv)}$, a translation.
Therefore, we have $T_{-\PtAv -\tau(\PtAv)}^{*}C = C$.
However, since $C$ is a principal polarization, we have $\Ker ~\varphi_{C} =\{0\}$.
Therefore, we have $T_{-\PtAv -\tau(\PtAv)} = id$. This means $\tau \circ T_{-\PtAv}$
is an automorphism of order 2 of $C$. By definition, this is not 
the inversion $\iota$. Hence, this gives an element of order 2 in ${\rm RA}(C)$.

(ii) $\Rightarrow$ (iii) This follows from Lemma \ref{replace}.

(iii) $\Rightarrow$ (i) This follows from Lemma \ref{freely} 
and Example \ref{Igusa}.
\qed}

\begin{remark} In the proof of the proposition, the automorphism 
$\tau \circ T_{-\PtAv}$ really gives a long element of order 2 in ${\rm RA}(C)$.
\end{remark}

By Castryck--Decru--Smith \cite[Subsection 3.3]{CDS}, 
if the curve $C$ of genus 2 is obtained
from a decomposed principally polarized abelian surface by a Richelot isogeny,
then the curve $C$ has a long reduced automorphism of order 2.
As is well-known, for a curve $C$ of genus 2, the Jacobian variety $J(C)$
has 15 Richelot isogenies (cf.~Castryck--Decru--Smith \cite[Subsection 3.2]{CDS}, 
for instance). 
If we have a Richelot isogeny $(A,C) \longrightarrow (A',C')$, then we also have a Richelot isogeny $(A',C') \longrightarrow (A,C)$. 
Therefore, we have the following proposition.

\begin{proposition}\label{decomposed}
Let $C$ be a nonsingular projective curve of genus 2.
Among the 15 Richelot isogenies outgoing from $J(C)$, the number of decomposed
Richelot isogenies is equal to the number of long elements
of order 2 in ${\rm RA}(C)$.
\end{proposition}

In this proposition, we consider that a different isotropic subgroup gives
a different Richelot isogeny. However, two different Richelot isogenies
may be isomorphic to each other by a suitable automorphism (see Definition \ref{Def:IsomorphismRichelot}).
From the next section, we will compute the number of Richelot isogenies up to
isomorphism.


\section{The number of long elements of order 2}
\label{Sec:NumberOfLongElements}

In this section, we count the number of long elements of order 2 
in ${\rm RA}(C)$.
For an element $f \in {\rm RA}(C)$, we express the reduced automorphism by
$$
    f : x \mapsto f(x)
$$
with a suitable coordinate $x$ of ${\bf A}^{1}\subset {\bf P}^{1}$.
We will give the list of $f(x)$ corresponding to elements of order 2. 
Here, we denote by $\omega$ a primitive cube root of unity, 
by $i$ a primitive fourth root of unity,
and by $\zeta$ a primitive sixth root of unity.

\noindent
Case (0) ${\rm RA}(C) \cong \{0\}$.

There exist no long elements of order 2.

\noindent
Case (1) ${\rm RA}(C) \cong {\bf Z}/2{\bf Z}$.

The curve $C$ is given by $y^2 = (x^2 - 1)(x^2 - a^2)(x^2 - b^2)$.

There exists only one long element of order 2 given by $f(x) = -x$.

\noindent
Case (2) ${\rm RA}(C) \cong  S_{3}$.

The curve $C$ is given by $y^2 = (x^3 - 1)(x^3 - a^3)$.

There exist three  long elements of order 2 given by
$f(x) = \frac{a}{x}$, $\frac{\omega a}{x}$, $\frac{\omega^2 a}{x}$.

\noindent
Case (3) ${\rm RA}(C) \cong  {\bf Z}/2{\bf Z}\times {\bf Z}/2{\bf Z}$.

The curve $C$ is given by $y^2 = x(x^2 - 1)(x^2 - a^2)$.

There exist  two long elements of order 2 given by
$f(x) = \frac{a}{x}$, $\frac{-a}{x}$,

and  there exists one short element of order 2 given by
$f(x) = -x$.

\noindent
Case (4) ${\rm RA}(C) \cong  D_{12}$.

The curve is given by $y^2 = x^6 - 1$.

There exist four long elements of order 2 given by
$f(x) = - x$, $\frac{\zeta}{x}$, $\frac{\zeta^3}{x}$, $\frac{\zeta^5}{x}$,

and there exist three short elements of order 2 given by
$f(x) = \frac{1}{x}$, $\frac{\zeta^2}{x}$, $\frac{\zeta^4}{x}$.

\noindent
Case (5) ${\rm RA}(C) \cong  S_{4}$.

The curve $C$ is given by $y^2 = x(x^4 - 1)$.

There exist six long elements of order 2 given by
$f(x)$$=$$\frac{x + 1}{x - 1}$,$- \frac{x - 1}{x + 1}$,$\frac{i(x + i)}{x - i}$,$\frac{i}{x}$,$- \frac{i}{x}$,$- \frac{i(x - i)}{x + i}$,

and there exist three short elements of order 2 given by
$f(x) =  - x$, $\frac{1}{x}$, $ - \frac{1}{x}$.

\noindent
Case (6) ${\rm RA}(C) \cong  {\bf Z}/5{\bf Z}$.

The curve is given by $y^2 = x^5 -1$.
 
There exist no long elements of order 2.
 
\begin{remark}
By Proposition \ref{decomposed} and the calculation above, we see
that for a curve $C$ of genus 2, the number of outgoing decomposed Richelot
isogenies from $J(C)$ is at most six. This result coincides with the one
given by Castryck--Decru--Smith \cite[Theorem 2]{CDS}.
\label{Rem:Thm3CDS}
\end{remark}


\section{Counting Richelot isogenies}
\label{Sec:Counting}


\subsection{Richelot isogenies from Jacobians of irreducible genus-2 curves}
\label{Sec:LociReducedAuto}

Let $C$ be a nonsingular projective curve of genus 2, and let 
$J(C)$ be the Jacobian variety of $C$. 
For a fixed $C$, we consider the set $\{(J(C), G)\}$
of pairs of $J(C)$ and an isotropic subgroup 
$G$ for the polarization $2C$. The group ${\rm Aut}(C)$ 
 acts on the ramification points of 
$C\longrightarrow {\bf P}^{1}$.
Using this action, ${\rm Aut}(C)$ induces the action 
on the set $\{(J(C), G)\}$. Since the inversion $\iota$ 
of $C$ acts on $J(C)[2]$ trivially,
the reduced group ${\rm RA}(C)$ of automorphisms acts on 
the set $\{(J(C), G)\}$ which consists of 15 elements.

Let $P_i$ $(i = 1, 2, \ldots, 6)$ be
the ramification points of $\psi : C \longrightarrow {\bf P}^{1}$. 
A division into the sets of 3 pairs
of these 6 points gives an isotropic subgroup $G$,
that is, 
$$
\{ P_{i_1} - P_{i_2}, P_{i_3} - P_{i_4}, P_{i_5} - P_{i_6}, {\rm the~ identity}\} 
$$
gives an isotropic subgroup of $J(C)[2]$. 
The action of ${\rm RA}(C)$ on the set $\{(J(C), G)\}$ is
given by the action of ${\rm RA}(C)$ on the set 
$$
\{\langle (P_{i_1}, P_{i_2}), (P_{i_3}, P_{i_4}), (P_{i_5}, P_{i_6})\rangle \},
$$
which contains 15 sets. Here, the pair $(P_{i}, P_{j})$ is unordered.
In this section, we count the number of orbits of this action
for each case. 

Let $C$ be a curve of genus 2 with ${\rm RA}(C) \cong {\bf Z}/2{\bf Z}$.
Such a curve is given by the equation 
$$
   y^2 = (x^2 - 1)(x^2 - a)(x^2 - b)
$$
with suitable conditions for $a$ and $b$.
The branch points $Q_i = \psi(P_i)$ are given by
$$
Q_1 = 1, ~Q_2 = -1, ~Q_3 = \sqrt{a}, ~Q_4 = - \sqrt{a}, ~Q_5 = \sqrt{b}, 
~ Q_6 = -\sqrt{b}.
$$
The generator of the reduced group ${\rm RA}(C)$ of automorphisms is given by
$$
\sigma : x \mapsto -x.
$$
Since the inversion $\iota$ acts trivially on the ramification points,
${\rm RA}(C)$ acts on the set of the ramification points
$\{ P_{1}, P_{2}, P_{3}, P_{4}, P_{5}, P_{6}\}$,
and the action of $\sigma$ on the ramification points is given by
$$ 
P_{2i- 1} \mapsto P_{2i},~ P_{2i} \mapsto P_{2i -1}\quad (i = 1, 2, 3).
$$
The isotropic subgroup which corresponds to 
$\langle (P_1, P_2), (P_3, P_4), (P_5, P_6)\rangle$ 
gives  a decomposed Richelot isogeny and the other isotropic subgroups
give  non-decomposed isogenies. Moreover, 
$\langle (\sigma(P_{i_1}), \sigma(P_{i_2})), (\sigma(P_{i_3}), \sigma(P_{i_4})), 
(\sigma(P_{i_5}), \sigma(P_{i_6}))\rangle$ gives the Richelot
isogeny isomorphic to the one given by
$\langle (P_{i_1}, P_{i_2}), (P_{i_3}, P_{i_4}), (P_{i_5}, P_{i_6})\rangle $.
We denote $P_i$ by $i$ for the sake of simplicity. Then, the action
$\sigma$ is given by the permutation $(1,2)(3,4)(5, 6)$, and by the action
of ${\rm RA}(C)$, the set 
$\{\langle (P_{i_1}, P_{i_2}), (P_{i_3}, P_{i_4}), (P_{i_5}, P_{i_6})\rangle \}$ 
of 15 elements is divided into the following 11 loci:
$$
\begin{array}{l}
\{[(1,2), (3, 4), (5, 6)]\},~
\{[(1,2), (3, 5), (4, 6)]\},~
\{[(1,2), (3, 6), (4, 5)]\},\\
\{[(1,3), (2, 4), (5, 6)]\},~
\{[(1,3), (2, 5), (4, 6)], [(1,6), (2, 4), (3, 5)]\}, \\
\{[(1,3), (2, 6), (4, 5)], [(1,5), (2, 4), (3, 6)]\},~
\{[(1,4), (2, 3), (5, 6)]\}, \\
\{[(1,4), (2, 5), (3, 6)], [(1,6), (2, 3), (4, 5)]\},
\{[(1,4), (2, 6), (3, 5)], [(1,5), (2, 3), (4, 6)]\}, \\
\{[(1,5), (2, 6), (3, 4)]\},~
\{[(1,6), (2, 5), (3, 4)]\}.
\end{array}
$$
The reduced automorphism $\sigma$ is a long one of order 2 
and the element $[(1,2), (3, 4), (5, 6)]$ is  
fixed by $\sigma$. Therefore, the element $[(1,2), (3, 4), (5, 6)]$
gives a decomposed isogeny. The other 10 loci give non-decomposed isogenies.
In the same way, we have the following proposition.

\begin{proposition}\label{irreducible}
Under the notation above, the number of Richelot isogenies
up to isomorphism in each case and the number of elements in each orbit 
are listed as follows. Here, in the list, for example,
$(1\times 6, 2 \times 4)(1 \times 1)$ means that there exist 6 orbits which 
contain 1 element and 4 orbits which contain 2 elements for non-decomposed
Richelot isogenies, and there exists 1 orbit which contains 1 element
for decomposed Richelot isogenies.

(0) ${\rm RA}(C) \cong \{0\}$: 15 Richelot isogenies. No decomposed one. 

\quad $(1 \times 15)(0)$.

(1) ${\rm RA}(C) \cong {\bf Z}/2{\bf Z}$: 11 Richelot isogenies. 1 decomposed one.

\quad $(1\times 6, 2 \times 4)(1 \times 1)$.

(2) ${\rm RA}(C) \cong S_3$: 7 Richelot isogenies. 1 decomposed one.

\quad $(1\times 3, 3 \times 3)(3 \times 1)$.

(3) ${\rm RA}(C) \cong {\bf Z}/2{\bf Z} \times {\bf Z}/2{\bf Z}$: 8 Richelot isogenies. 2 decomposed ones. 

\quad $(1\times 1, 2 \times 4, 4\times 1)(1 \times 2)$.

(4) ${\rm RA}(C) \cong D_{12}$: 5 Richelot isogenies. 2 decomposed ones. 

\quad $(2\times 1, 3 \times 1, 6 \times 1)(1 \times 1, 3\times 1)$.

(5) ${\rm RA}(C) \cong S_4$: 4 Richelot isogenies. 1 decomposed one. 

\quad $(1\times 1, 4 \times 2)(6 \times 1)$.

(6) ${\rm RA}(C) \cong {\bf Z}/5{\bf Z}$: 3 Richelot isogenies. No decomposed one.

\quad $(5\times 3)(0)$.

\end{proposition}


\begin{theorem}
\label{Thm:RichelotFromNonDecomp}
The total number of Richelot isogenies up to isomorphism outgoing 
from the irreducible superspecial curves of genus 2 is equal to
$$
\frac{(p - 1)(p + 2)(p + 7)}{192} - 3\{1 - (\frac{-1}{p})\}/32 + \{1 - (\frac{-2}{p})\}/8.
$$
The total number of decomposed Richelot isogenies up to isomorphism 
outgoing from the irreducible superspecial curves of genus 2 is equal to
\begin{equation}
\label{Eq:NumberOfDecompRichelot}
\frac{(p-1)(p+3)}{48} - \{1 - (\frac{-1}{p})\}/8 + \{1 - (\frac{-3}{p})\}/6.
\end{equation}
\end{theorem}

\proof{The total number of Richelot isogenies up to isomorphism 
outgoing 
from the irreducible superspecial curves of genus 2 is equal to
$$
15n_0 + 11n_1 + 7n_2 + 8n_3 + 5n_4 + 4n_5 + 3n_6
$$
and the total number of decomposed Richelot isogenies up to isomorphism
outgoing from the irreducible superspecial curves of genus 2 is 
equal to
$$
n_1 + n_2 + 2n_3 + 2n_4 + n_5.
$$
The results follow from these facts.
\qed}


\subsection{Richelot isogenies from elliptic curve products}
\label{Sec:DecomposedCases}

Let $E$, $E'$ be supersingular elliptic curves, and we consider
a decomposed principal polarization $E + E'$ and
a Richelot isogeny $(E \times E', E + E') \longrightarrow (J(C), C)$.
For a principally polarized abelian surface $(E \times E', E + E')$,
we denote by ${\rm Aut}(E\times E')$ the group of automorphisms of
$E\times E'$ which preserve the polarization $E + E'$. 
Let $\{P_1, P_2, P_3\}$ (resp.\,$\{P_4, P_5, P_6\}$) be the 2-torsion points of
$E'$ (resp.\,$E$). Then, the six points $P_i$ ($1 \leq i \leq 6$) on $E \times E'$
play the role of ramification points of irreducible curves of genus 2, and
${\rm Aut}(E\times E')$ acts on the set $\{P_1, P_2, P_3, P_4, P_5, P_6\}$. 
Note that the subgroup 
$\langle \iota_E\times id_{E'}, id_{E}\times \iota_{E'} \rangle$ acts on
the set $\{P_1, P_2, P_3, P_4, P_5, P_6\}$ trivially.
In this section, let $E_2$ be the elliptic curve defined by
$
y^2 = x^3 - x
$
and $E_3$ the elliptic curve defined by
$
y^2 = x^3 - 1.
$
We know ${\rm Aut}{E_2}\cong {\bf Z}/4{\bf Z}$ and
${\rm Aut}{E_3}\cong {\bf Z}/6{\bf Z}$. The elliptic curve
$E_2$ is supersingular if and only if $p \equiv 3~({\rm mod}~4)$
and $E_3$ is supersingular if and only if $p \equiv 2~({\rm mod}~3)$.
In this section, the abelian surface $E \times E'$ means
an abelian surface $E \times E'$ with principal polarization $E + E'$.

Now, let $E$, $E'$ be supersingular elliptic curves
which are neither isomorphic to $E_2$ nor to $E_3$. We also assume 
$E$ is not isomorphic to $E'$. Using these notations, we have 
the following list of orders of the groups of automorphisms.
$$
\begin{array}{l}
|{\rm Aut}(E\times E')| = 4,~ 
|{\rm Aut}(E\times E)= 8,~ 
|{\rm Aut}(E\times E_2)|= 8, ~
 |{\rm Aut}(E\times E_3)|= 12,\\
  |{\rm Aut}(E_2\times E_2)| = 32, ~
 |{\rm Aut}(E_3\times E_3)| = 72,~ 
  |{\rm Aut}(E_2\times E_3)| = 24.
\end{array}
$$
The isotropic subgroups for the polarization $2(E + E')$ are
determined in Castryck--Decru--Smith \cite[Subsection 3.3]{CDS}.
Using their results and the same method as in Subsection \ref{Sec:LociReducedAuto}, we have 
the following proposition.

\begin{proposition}\label{reducible} Let $E$, $E'$ be supersingular 
elliptic curves
which are neither isomorphic to $E_2$ nor to $E_3$ with $E_2$ and $E_3$
defined as above. We also assume  that
$E$ is not isomorphic to $E'$. 
The number of Richelot isogenies up to isomorphism 
outgoing from a decomposed principally polarized
superspecial abelian surface in each case and the number of elements 
in each orbit are listed as follows. 
Here, in the list, for example,
$(1\times 3, 2 \times 1)(1 \times 4, 2\times 3)$ means that 
there exist 3 orbits which 
contain 1 element and 1 orbit which contains 2 elements for non-decomposed
Richelot isogenies, and there exist 4 orbits which contain 1 element
and 3 orbits which contain 2 elements for decomposed Richelot isogenies.

(i) $E\times E'$ : 15 Richelot isogenies, 6 non-decomposed ones. 

\quad $(1\times 6)(1\times 9)$.

(ii) $E\times E$ : 11 Richelot isogenies, 4 non-decomposed ones.

\quad $(1\times 3, 2\times 1)(1\times 4, 2\times 3)$.

(iii) $E \times E_2$ : 9 Richelot isogenies, 3 non-decomposed ones $(p \equiv 3~({\rm mod}~4))$. 

\quad $(2\times 3)(1\times 3, 2\times 3)$.

(iv) $E\times E_3$ : 5 Richelot isogenies, 2 non-decomposed ones $(p \equiv 2~({\rm mod}~3))$. 

\quad $(3\times 2)(3\times 3)$.

(v) $E_2\times E_2$ : 5 Richelot isogenies, 1 non-decomposed one $(p \equiv 3~({\rm mod}~4))$. 

\quad $(4\times 1)(1\times 1, 2\times 1, 4\times 2)$.

(vi) $E_3\times E_3$ : 3 Richelot isogenies, 1 non-decomposed one $(p \equiv 2~({\rm mod}~3))$. 

\quad $(3\times 1)(3\times 1, 9\times 1)$.

(vii) $E_2\times E_3$ : 3 Richelot isogenies, 1 non-decomposed one $(p \equiv 11~({\rm mod}~12))$. 

\quad $(6\times 1)(3\times 1, 6\times 1)$.
\end{proposition}
\proof{
We give a proof for the case (iv). 
For the other cases, the arguments are quite similar.
Since the elliptic curve $E_3$ is defined by $y^2 = x^3 - 1$,
the 2-torsion points $(x, y)$ of $E_3$ are given by
$P_1 = (1, 0)$, $P_2 = (\omega, 0)$ and $P_3 =(\omega^2, 0)$.
Here, $\omega$ is a primitive cube root of unity. 
We denote by $P_4$, $P_5$ and $P_6$ the 2-torsion points of $E$.
We have an automorphism $\sigma$ of order 3 of $E_3$ defined by
$
   \sigma : x \mapsto \omega x, y \mapsto y.
$ 
As in the case of Subsection 6.1, we describe the isotropic subgroups $G$.
We know that a division into the sets of 3 pairs
of these 6 points $P_{i}$ ($1 \leq i \leq 6$) on $E \times E_3$ gives an isotropic subgroup $G$,
that is, 
$
\{ P_{i_1} - P_{i_2}, P_{i_3} - P_{i_4}, P_{i_5} - P_{i_6}, {\rm the~ identity}\} 
$
gives an isotropic subgroup of $(E\times E_3)[2]$. Here, we consider 
$P_i$ ($1 \leq i \leq 3$) as the point $(0, P_i)$ on $E \times E_3$, and 
$P_i$ ($4 \leq i \leq 6$) as the point $(P_i, 0)$ on $E \times E_3$.
This set contains 
15 elements.
In the case (iv), we have $E \not\cong E_3$. Therefore,
by Castryck--Decru--Smith \cite[Subsection 3.3]{CDS}, 
among the 15 isotropic subgroups 
the 9 cases such that $P_{i_1}, P_{i_2}, P_{i_3} \in E$ and 
$P_{i_4}, P_{i_5}, P_{i_6} \in E_3$ give the decomposed Richelot isogenies
and the rest gives the non-decomposed Richelot isogenies.
For the abbreviation, we denote by $P_i$ by $i$. Then, on the set
$\{1, 2, 3, 4, 5, 6\}$, $id_E \times \sigma$ acts 
as the cyclic permutation $(1, 2, 3)$.
The isotropic subgroup $G$ is determined by the set
of 3 pairs of 2-torsion points:
$$
\{({i_1}, {i_2}), ({i_3}, {i_4}), ({i_5}, {i_6})\},
$$
and the group ${\rm Aut}(E \times E_3)$
induces the action on the set of the 15 isotropic subgroups. 
Since the action of the subgroup $\langle \iota_E\times id_{E_3}, id_{E}\times \iota_{E_3} \rangle$ is trivial on the set 
of the 15 isotropic subgroups, we see that the action is
given by the group 
${\rm Aut}(E \times E_3)/\langle \iota_E\times id_{E_3}, id_{E}\times \iota_{E_3} \rangle\cong \langle id_E\times \sigma\rangle$. 
By this action, the set of the 15 isotropic subgroups is divided 
into the following 5 orbits:
$$
\begin{array}{l}
\{[(1, 2), (3, 4), (5, 6)], [(2, 3), (1, 4), (5, 6)], [(1, 3), (2, 4), (5, 6)]\},\\
\{[(1, 2), (3, 5), (4, 6)], [(2, 3), (1, 5), (4, 6)], [(1, 3), (2, 5), (4, 6)]\},\\
\{[(1, 2), (3, 6), (4, 5)], [(2, 3), (1, 6), (4, 5)], [(1, 3), (2, 6), (4, 5)]\},\\
\{[(1, 4), (2, 5), (3, 6)], [(1, 6), (2, 4), (3, 5)], [(1, 5), (2, 6), (3, 5)]\},\\
\{[(1, 4), (2, 6), (3, 5)], [(1, 5), (2, 4), (3, 6)], [(1, 6), (2, 5), (3, 4)]\}.\\
\end{array}
$$
By the criterion above, the first 3 sets correspond with the decomposed Richelot
isogenies, and the last 2 sets correspond with the non-decomposed Richelot
isogenies. 
\qed}
\vspace*{0.3cm}

We denote by $h$ the number of supersingular elliptic curves defined over $k$.
Then, we know 
$$
h = \frac{p-1}{12} + \{1 - (\frac{-3}{p})\}/3 + \{1 - (\frac{-1}{p})\}/4
$$
(cf.~Igusa \cite{I0}, for instance).
We denote by $h_1$ the number of supersingular elliptic curves $E$
with ${\rm Aut}(E) \cong {\bf Z}/2{\bf Z}$, 
$h_2$ the number 
of supersingular elliptic curves $E_2$
with ${\rm Aut}(E_2) \cong {\bf Z}/4{\bf Z}$, 
$h_3$ the number of supersingular elliptic curves $E_3$
with ${\rm Aut}(E_3) \cong {\bf Z}/6{\bf Z}$.
We have $h = h_1 + h_2 + h_3$ and 
$
 h_2 = \{1 - (\frac{-1}{p})\}/2
$
and 
$
h_3 = \{1 - (\frac{-3}{p})\}/2
$.


\begin{theorem}
\label{Thm:RichelotFromDecomp}
The total number of non-decomposed Richelot isogenies up to isomorphism outgoing
from decomposed principally polorized superspecial abelian surfaces is equal to
\begin{equation}
\frac{(p-1)(p+3)}{48} - \{1 - (\frac{-1}{p})\}/8 + \{1 - (\frac{-3}{p})\}/6.
\label{Eq:NumberOfNonDecompRichelot}
\end{equation}
The total number of decomposed Richelot isogenies up to isomorphism outgoing
from decomposed principally polorized superspecial abelian surfaces is equal to
$$
\frac{(p - 1)(3p + 17)}{96} + (p + 6)\{1 - (\frac{-1}{p})\}/16 
+ \{1 - (\frac{-3}{p})\}/3.
$$
\end{theorem}

\proof{
The total number of non-decomposed Richelot isogenies up to isomorphism 
outgoing from decomposed principally polorized superspecial abelian surfaces 
is equal to
$$
6\{\frac{h_1(h_1 - 1)}{2}\} + 4 h_1 + 3h_2h_1 + 2h_3 h_1 + h_2 + h_3 + h_2h_3.
$$
The total number of decomposed Richelot isogenies up to isomorphism 
outgoing
from decomposed principally polorized superspecial abelian surfaces 
is equal to
$$
9\{\frac{h_1(h_1 - 1)}{2}\} + 7 h_1 + 6h_2h_1 + 3h_3 h_1 + 4h_2 + 2h_3 + 2h_2h_3.
$$
Since $\{1 - (\frac{-1}{p})\}^2 = 2\{1 - (\frac{-1}{p})\}$ and
$\{1 - (\frac{-3}{p})\}^2 = 2\{1 - (\frac{-3}{p})\}$, the result follows
from these facts.
\qed}

\begin{remark}
Since the total number of {\em decomposed} Richelot isogenies up to isomorphism 
outgoing from the {\em irreducible} superspecial curves of genus 2 is equal to 
the total number of {\em non-decomposed} Richelot isogenies up to isomorphism outgoing from {\em decomposed} principally polorized superspecial abelian surfaces, (\ref{Eq:NumberOfDecompRichelot}) and (\ref{Eq:NumberOfNonDecompRichelot}) give the same number.
\end{remark}


\section{Examples}
\label{Sec:Examples}

By Ibukiyama--Katsura--Oort \cite[Subsection 1.3]{IKO}, 
we have the following normal forms
of curves $C$ of genus 2 with given reduced group ${\rm RA}(C)$ of automorphims:

\begin{enumerate}
\item
For $S_3 \subset {\rm RA}(C)$, the normal form is
$
y^2 =(x^3 -1)(x^3 - \alpha).
$
This curve is superspecial if and only if $\alpha$ is a zero of
the polynomial
$$
g(\hVar) = \sum_{l = 0}^{[p/3]}
\left(
\begin{array}{c}
(p -1)/2 \\
((p + 1)/6) + l
\end{array}
\right)
\left(
\begin{array}{c}
(p -1)/2\\
 l
\end{array}
\right)
\hVar^l.
$$

\item
For ${\bf Z}/2{\bf Z} \times {\bf Z}/2{\bf Z} \subset {\rm RA}(C)$, the normal form is
$
y^2 = x(x^2 -1)(x^2 - \beta).
$
This curve is superspecial if and only if $\beta$ is a zero of
the polynomial
$$
h(\hVar) = \sum_{l = 0}^{[p/4]}
\left(
\begin{array}{c}
(p -1)/2 \\
((p + 1)/4) + l
\end{array}
\right)
\left(
\begin{array}{c}
(p -1)/2\\
 l
\end{array}
\right)
\hVar^l.
$$

\item
For ${\rm RA}(C) \cong D_{12}$, the normal form is
$
y^2 = x^6 -1.
$
This curve is superspecial if and only if $p \equiv 5 ~({\rm mod}~6)$
(cf.~Ibukiyama--Katsura--Oort \cite[Proposition 1.11]{IKO}).

\item
For ${\rm RA}(C) \cong S_4$, the normal form is
$
y^2 = x(x^4 -1).
$
This 
is superspecial if and only if $p \equiv 5 ~{\rm or}~7~({\rm mod}~8)$
(cf.~Ibukiyama--Katsura--Oort \cite[Proposition 1.12]{IKO}).
\end{enumerate}

Finally, the elliptic curve $E$ defined by
$
   y^2 = x(x -1)(x -\lambda)
$
is supersingular if and only if $\lambda$ is a zero of
the Legendre polynominal
$$
\Phi (\hVar) = \sum_{l = 0}^{(p -1)/2}
\left(
\begin{array}{c}
(p -1)/2 \\
    l
\end{array}
\right)^2
\hVar^l.
$$
Using these results, we construct some examples. 


\subsection{Examples in characteristic 13}
Assume the characteristic $p = 13$. 
Over $k$ we have
 only one supersingular elliptic curve $E$, and three superspecial
curves $C_1$, $C_2$ and $C_3$ of 
genus 2 with  
${\rm RA}(C_1) \cong S_3$, 
${\rm RA}(C_2) \cong {\bf Z}/2{\bf Z} \times {\bf Z}/2{\bf Z}$ 
and 
${\rm RA}(C_3) = S_4$, respectively 
(cf.~Ibukiyama--Katsura--Oort \cite[Remark 3.4]{IKO}). 
In characteristic 13, we know
$
h(\hVar) = 7\hVar^3 + 12\hVar^2 + 12\hVar + 7,
$
and the zeros are $-1$ and $-5 \pm \sqrt{6}$. 
We also know
$
g(\hVar) = 2\hVar^4 + 3\hVar^3 + 4\hVar^2 + 3\hVar + 2,
$
and one of the zeros is $-4 + \sqrt{2}$.
The Legendre polynomial is given by
$
\Phi(\hVar) = \hVar^6 + 10\hVar^5 + 4\hVar^4 + 10\hVar^3 + 4\hVar^2 + 10\hVar + 1,
$
and one of the zeros is  $3 - 2\sqrt{2}$. Using these facts, we know
that the curves above are given by the following equations:

(1) $E$: $y^2 = x(x - 1)(x - 3 + 2\sqrt{2})$ \quad $({\rm RA}(E) ={\rm Aut}(E)/\langle\iota_{E}\rangle \cong \{0\}$),
 
(2) $C_1$: $y^2 = (x^3 - 1)(x^3 + 4 -\sqrt{2})$ \quad $({\rm RA}(C_1) \cong S_3)$, 

(3) $C_2$: $y^2 = x(x^2 - 1)(x^2 + 5  + 2\sqrt{6})$ \quad 
$({\rm RA}(C_2) \cong {\bf Z}/2{\bf Z} \times {\bf Z}/2{\bf Z})$,

(4) $C_3$: $y^2 = x(x^4 - 1)$ \hspace{0.2cm}$({\rm RA}(C_3) \cong S_4)$.

\vspace{-0.5cm}
\begin{figure}[hbt]
\hspace*{-0.9cm}
\begin{center}
\includegraphics[scale=0.40, clip]{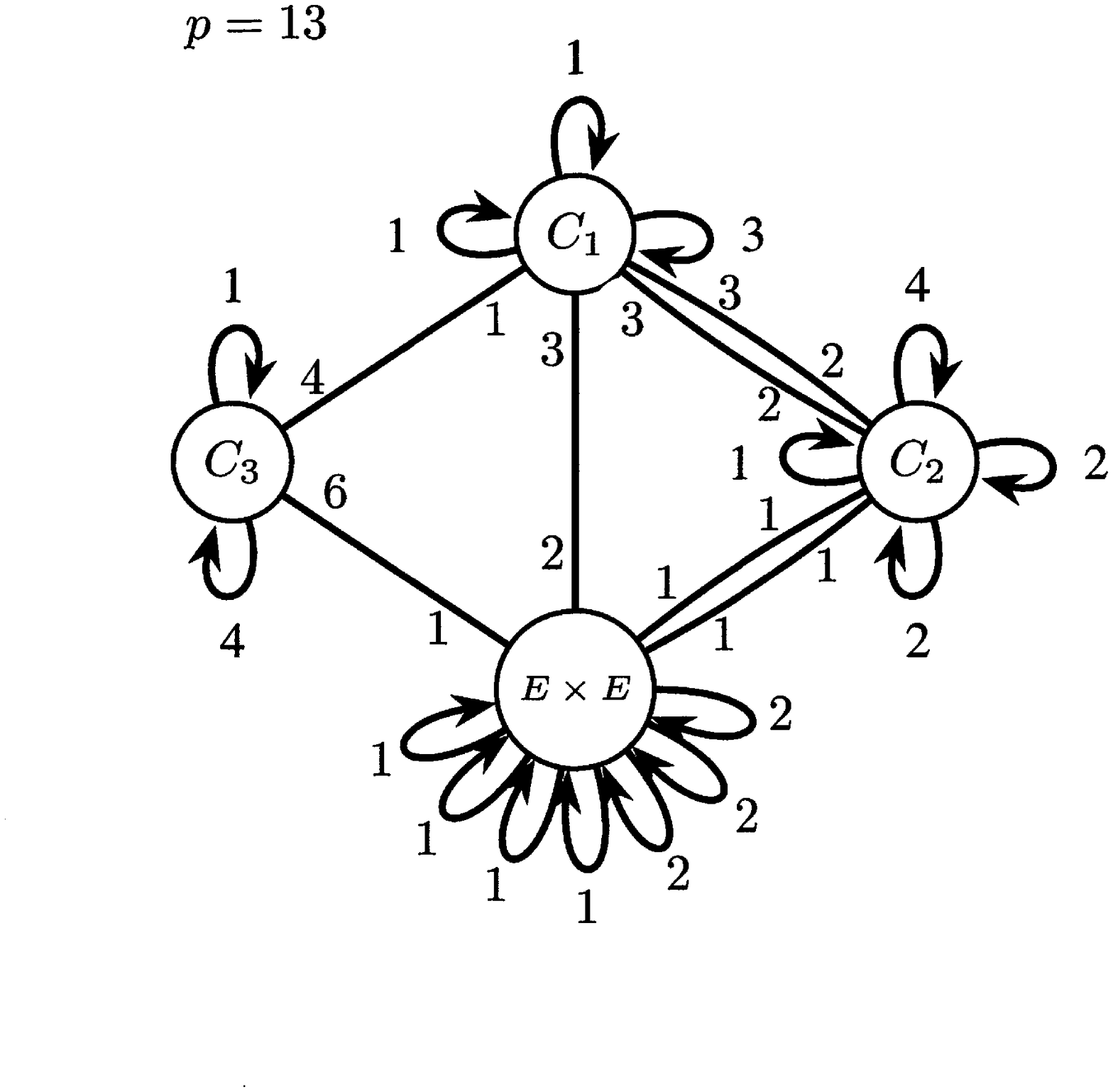}
\end{center}
\end{figure}

\vspace{-1.5cm}
\noindent
Therefore, outgoing from superspecial curves of genus 2,
we have, in total, $1 + 2 + 1 = 4$ decomposed
Richelot isogenies up to isomorphism by Proposition \ref{irreducible}. 
On the other hand, outgoing from the unique 
decomposed principally polarized
abelian surface $(E \times E, E + E)$, we have $5$ non-decomposed Richelot
isogenies (not up to isomorphism) 
(cf.~Igusa \cite{I0} and Castryck--Decru--Smith \cite[Figure 1]{CDS}). 
Using the method in 
Castryck--Decru--Smith \cite[Subsection 3.3]{CDS},
as the images of $5$ non-decomposed Richelot isogenies, we have the following superspecial curves
of genus 2:

(a) $C_a$: $y^2 = (x^2 - 1)(x^2 - 4 + 7\sqrt{2})(x^2 - 6 + 6\sqrt{2})$ \quad  
$({\rm RA}{(C_a)} \cong {\bf Z}/2{\bf Z} \times {\bf Z}/2{\bf Z})$,

(b) $C_b$: $y^2 = (x^2 - 1)(x^2 + 3 - 2\sqrt{2})(x^2 - 4 - \sqrt{2})$ \quad  
$({\rm RA}{(C_b)}  \cong S_4)$,

(c) $C_c$: $y^2 = (x^2 - 1)(x^2 + 3 - 4\sqrt{2})(x^2 + 1 + 3\sqrt{2})$ \quad 
$({\rm RA}{(C_c)}  \cong S_3)$,

(d) $C_d$: $y^2 = (x^2 - 1)(x^2 - 3)(x^2 + 3 - 4\sqrt{2})$  \quad 
$({\rm RA}{(C_d)}  \cong S_3)$,

(e) $C_e$: $y^2 = (x^2 - 1)(x^2 - 6 - 6\sqrt{2})(x^2 - 2 + 2\sqrt{2})$ \quad  
$({\rm RA}{(C_e)} \cong {\bf Z}/2{\bf Z} \times {\bf Z}/2{\bf Z})$.

\noindent
We see that $C_a \cong C_e \cong C_2$, $C_c \cong C_d \cong C_1$ and $C_b \cong C_3$.
As Richelot isogenies, $(E \times E, E + E) \longrightarrow (J(C_c), C_c)$ is isomorphic 
to $(E \times E, E + E) \longrightarrow (J(C_d), C_d)$, but 
$(E \times E, E + E) \longrightarrow (J(C_a), C_a)$ is not isomorphic 
to $(E \times E, E + E) \longrightarrow (J(C_e), C_e)$. 
Compare our graph with Castryck--Decru--Smith \cite[Figure 1]{CDS}.
In the graph the numbers along the edges are the multiplicities of 
Richelot isogenies outgoing from the nodes.


\subsection{Examples in characteristic 11}\label{Sec:Examples11}
Assume the characteristic $p = 11$. Over $k$ we have
two supersingular elliptic curves $E_2$,$E_3$ and 
two superspecial
curves $C_1$, $C_2$ of genus 2 with 
${\rm RA}(C_1) \cong S_3$, ${\rm RA}(C_2) \cong D_{12}$, respectively
(cf.~Ibukiyama--Katsura--Oort \cite[Remark 3.4]{IKO}). 
In characteristic 11, we know
$$
g(\hVar) = 10(\hVar^3 + 5\hVar^2 + 5\hVar + 1),
$$
and the roots are $-1$, $3$ and $4$. Using this fact, we know
that the curves

\noindent
above are given by the following equations:

(1) $E_2$: $y^2 = x^3-x$ \quad
$({\rm RA}(E_2) \cong {\bf Z}/2{\bf Z}$),

(2) $E_3$: $y^2 = x^3 - 1$ \quad 
$({\rm RA}(E_3) \cong {\bf Z}/3{\bf Z})$, 

(3) $C_1$: $y^2 = (x^3 - 1)(x^3 -3)$ \quad 
$({\rm RA}(C_1) \cong S_3)$,

(4) $C_2$: $y^2 = x^6 - 1$ \quad $({\rm RA}(C_2) \cong D_{12})$.

\vspace{-0.5cm}
\begin{figure}[hbt]
\hspace*{-0.9cm}
\begin{center}
\includegraphics[scale=0.40,clip]{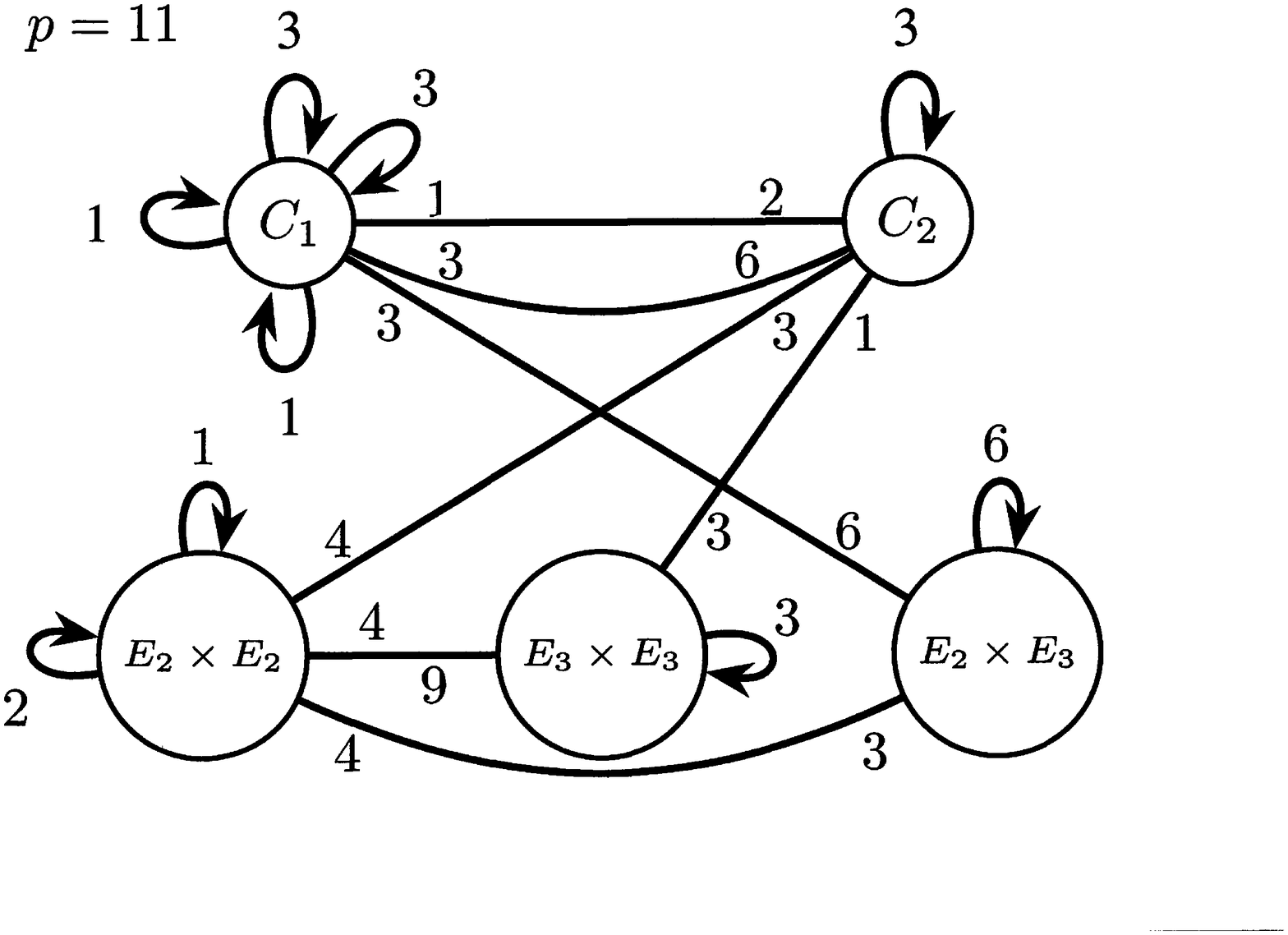}
\end{center}
\end{figure}

\vspace{-1.2cm}
\noindent
We have three decomposed principally polarized abelian surfaces:
$
E_2\times E_2, E_3\times E_3, E_2 \times E_3.
$
Therefore, from the superspecial curves of genus 2 we have, in total, $1 + 2 = 3$ decomposed
Richelot isogenies up to isomorphism by Proposition \ref{irreducible}. 
On the other hand, from the decomposed principally polarized
abelian surfaces, we have $1 + 1 + 1 = 3$ non-decomposed Richelot
isogenies up to isomorphism by Proposition \ref{reducible}.
For the decomposed 
principally polarized abelian surface $E_2\times E_2$ 
the image of the only one non-decomposed Richelot isogeny is given by $C_2$.
For the decomposed principally polarized abelian surface $E_3\times E_3$ 
the image of the only one non-decomposed Richelot isogeny is also given by $C_2$.
For the decomposed principally polarized abelian surface $E_2\times E_3$ 
the image of the only one non-decomposed Richelot isogeny is given by $C_1$.
See also Jordan--Zaytman \cite[Subsection 5.1]{JZ}.


\subsection{Examples in characteristic 7}
Assume the characteristic $p = 7$. 
Over $k$ we have  only one supersingular 
elliptic curve $E_2$ and 
only one superspecial
curves $C$ of genus 2, which has  
${\rm RA}(C) \cong S_4$
(cf.~Ibukiyama--Katsura--Oort \cite[Remark 3.4]{IKO}). 
They are given by the following equations:

\quad (1) $E_2$: $y^2 = x^3-x$ \quad $({\rm RA}(E_2) \cong {\bf Z}/2{\bf Z}$),

\quad (2) $C$: $y^2 = x(x^4 - 1)$ \quad $({\rm RA}(C) \cong S_4)$.

We have only one decomposed principally polarized abelian surface
$
E_2\times E_2.
$
Therefore, outgoing from the superspecial curves of genus 2 
we have only one decomposed
Richelot isogeny up to isomorphism. From the decomposed principally 
polarized abelian surface, we also have only one non-decomposed 
Richelot isogeny up to isomorphism 
(cf.~Castryck--Decru--Smith \cite[Subsections 3.2 and 3.3]{CDS}). 
For the decomposed principally polarized abelian surface $E_2\times E_2$ 
the image of the only one non-decomposed Richelot isogeny is given by $C$.

\vspace{-0.5cm}
\begin{figure}[hbt]
\hspace*{-0.9cm}
\begin{center}
\includegraphics[scale=0.40,clip]{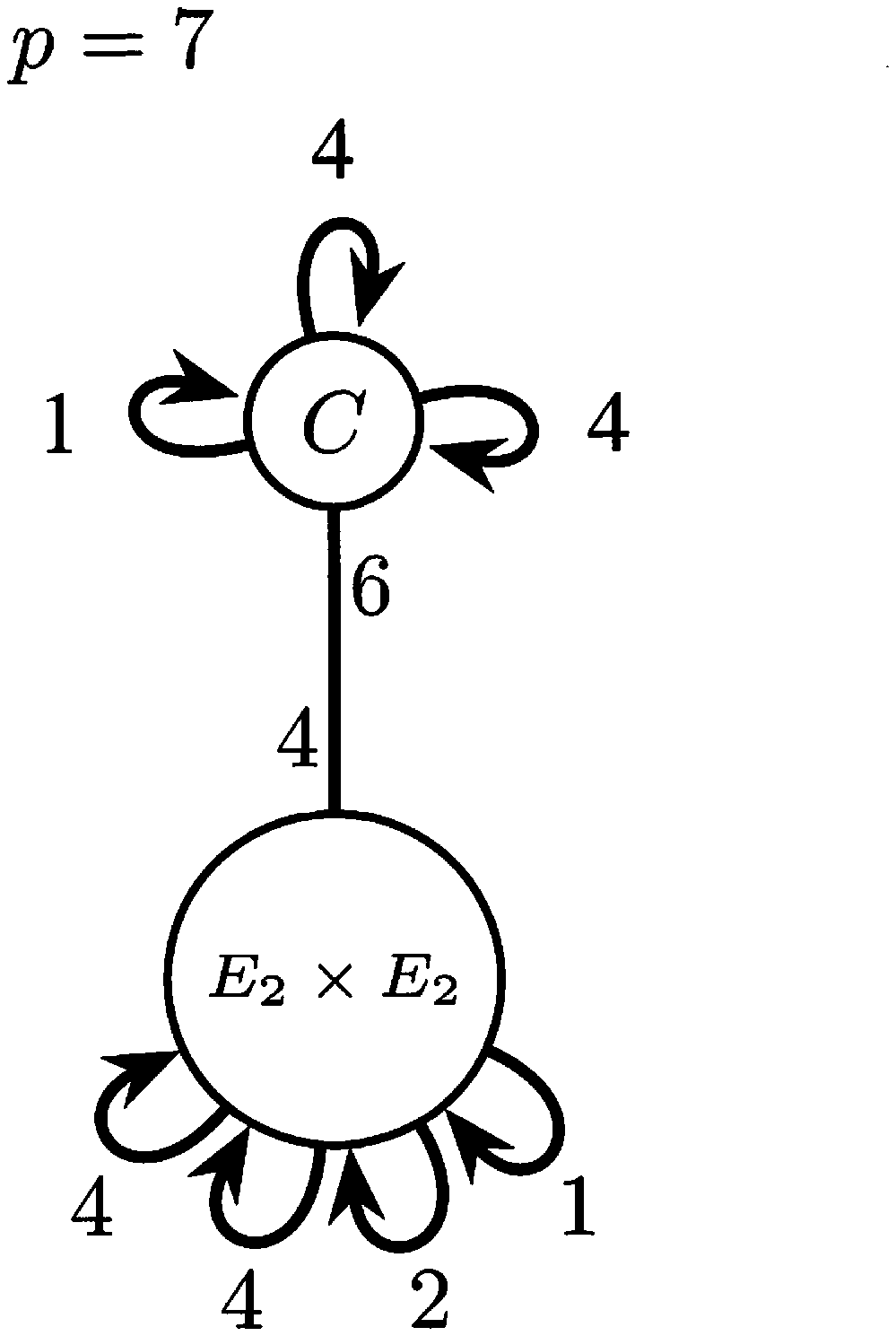}
\end{center}
\end{figure}
\vspace{-1.5cm}

\section{Concluding Remark}
\label{Sec:ConcludingRemarks}
Our results answered a question about the number of decomposed Richelot isogenies and 
improved our understanding of the isogeny graph for genus-2 isogeny cryptography. Further applications (or implications) of our results to cryptography are left as an open problem. 

For example, a very recent cryptanalytic algorithm by Costello and Smith \cite{CS20} is considered as an interesting target. 
They reduced the isogeny path-finding algorithm in the superspecial Richelot isogeny graph to the elliptic curve path-finding problem, thus improving the complexity. A key ingredient of the reduction is a sub-algorithm for finding a path connecting a given irreducible genus-2 curve and the (connected) subgraph consisting of elliptic curve products.

Proposition \ref{main} showed the equivalence of existence of a decomposed Richelot isogeny outgoing from $J(C)$ and 
that of a (long) element of order 2 in the reduced group of automorphisms of $C$. It implies that the 
subgraph of elliptic curve products are adjacent to genus-2 curves having involutive reduced automorphisms in the superspecial graph. We hope that this new characterization can be applied to analysing and/or improving the Costello--Smith attack.  

\bibliographystyle{abbrv}
\bibliography{ANTS2020}

\begin{thebibliography}{10}

\bibitem{CDS}
W.~Castryck, T.~Decru, and B.~Smith.
\newblock Hash functions from superspecial genus-2 curves using {R}ichelot
  isogenies.
\newblock In {\em {NutMiC 2019: Number-Theoretic Methods in Cryptology}}, 2019.
\newblock To appear in J.\,of Math.\,Crypt.

\bibitem{CLMPR18}
W.~Castryck, T.~Lange, C.~Martindale, L.~Panny, and J.~Renes.
\newblock {CSIDH}: An efficient post-quantum commutative group action.
\newblock In {\em {ASIACRYPT} 2018, Part {III}}, pages 395--427, 2018.

\bibitem{CLG09}
D.~Charles, K.~Lauter, and E.~Goren.
\newblock Cryptographic hash functions from expander graphs.
\newblock {\em J.~Crypt.}, 22(1):93--113, 2009.

\bibitem{CS20}
C.~Costello and B.~Smith.
\newblock The supersingular isogeny problem in genus 2 and beyond.
\newblock In {\em PQCrypto 2020}, pages 151--168, 2020.

\bibitem{FJP14}
L.~D. Feo, D.~Jao, and J.~Pl{\^{u}}t.
\newblock Towards quantum-resistant cryptosystems from supersingular elliptic
  curve isogenies.
\newblock {\em J.~Math.~Crypt.}, 8(3):209--247, 2014.

\bibitem{FT19}
E.~V. Flynn and Y.~B. Ti.
\newblock Genus two isogeny cryptography.
\newblock In {\em PQCrypto 2019}, pages 286--306, 2019.

\bibitem{IKO}
T.~Ibukiyama, T.~Katsura, and F.~Oort.
\newblock {Supersingular curves of genus two and class numbers}.
\newblock {\em Compositio Math.}, 57:127--152, 1986.

\bibitem{I0}
J.-I. Igusa.
\newblock {Class number of a definite quaternion with prime discriminant}.
\newblock {\em Proc. Nat. Acad. Sci. U.S.A.}, 44:312--314, 1958.

\bibitem{I}
J.-I. Igusa.
\newblock {Arithmetic variety of moduli for genus two}.
\newblock {\em Ann of Math.}, 72:612--649, 1960.

\bibitem{SIKE20}
D.~Jao, R.~Azarderakhsh, M.~Campagna, C.~Costello, L.~D. Feo, B.~Hess,
  A.~Jalali, B.~Koziel, B.~La{M}acchia, P.~Longa, M.~Naehrig, G.~Pereira,
  J.~Renes, V.~Soukharev, and D.~Urbanik.
\newblock {SIKE}: Supersingular isogeny key encapsulation.
\newblock {\em submission to the NIST's PQC standardization, round 2, updated
  version}, April 2020.

\bibitem{JZ}
B.~W. Jordan and Y.~Zaytman.
\newblock Isogeny graphs of superspecial abelian varieties and generalized
  {B}randt matrices.
\newblock {\em ArXiv}, abs/2005.09031, 2020.

\bibitem{KO2}
T.~Katsura and F.~Oort.
\newblock {Families of supersingular abelian surfaces}.
\newblock {\em Compositio Math.}, 62:107--167, 1987.

\bibitem{KO}
T.~Katsura and F.~Oort.
\newblock {Supersingular abelian varieties of dimension two or three and class
  numbers}.
\newblock {\em Advanced Studies in Pure Math.}, 10:253--281, 1987.

\bibitem{M}
D.~Mumford.
\newblock {\em Abelian Varieties}.
\newblock Oxford Univ. Press, 1970.

\bibitem{S}
T.~Shioda.
\newblock {Supersingular K3 surfaces}.
\newblock In {\em Algebraic Geometry, Proc. Copenhagen 1978 (K. L$\phi$nsted,
  ed.)}, Lecture Notes in Math. 732, pages 563--591. Springer-Verlag,
  Berlin-Heidelberg-New York, 1979.

\bibitem{Sm}
B.~Smith.
\newblock {\em Explicit endomorphisms and correspondences}.
\newblock PhD thesis, The University of Sydney, 2005.

\bibitem{Tak17}
K.~Takashima.
\newblock Efficient algorithms for isogeny sequences and their cryptographic
  applications.
\newblock In {\em Mathematical Modelling for Next-Generation Cryptography:
  {CREST} Crypto-Math Project}, pages 97--114. {Springer Verlag}, 2017.

\end{thebibliography}

\end{document}